\documentclass[12pt]{article}
\usepackage{multicol}
\usepackage{amssymb}
\usepackage{amsfonts}
\usepackage{amsmath}
\usepackage{amsbsy}
\usepackage[dvips]{graphics}
\usepackage{graphics,pspicture}
\usepackage[small,nohug,heads=littlevee]{diagrams}
\diagramstyle[labelstyle=\scriptstyle]
\newcommand{\be}{\begin{equation}}
\newcommand{\ee}{\end{equation}}
\newcommand{\bea}{\begin{eqnarray}}
\newcommand{\eea}{\end{eqnarray}}
\newcommand{\ba}{\begin{array}}
\newcommand{\ea}{\end{array}}

\newcommand{\bc}{\begin{center}}
\newcommand{\ec}{\end{center}}
\newcommand{\ben}{\begin{enumerate}}
\newcommand{\een}{\end{enumerate}}
\newcommand{\bfi}{\begin{figure}}
\newcommand{\efi}{\end{figure}}

\newcommand{\bq}{\begin{quote}}
\newcommand{\eq}{\end{quote}}
\newcommand{\bqu}{\begin{quotation}}
\newcommand{\equ}{\end{quotation}}
\newenvironment{emphit}{\begin{itemize}}{\end{itemize}}
\newcommand{\bemp}{\begin{emphit}}
\newcommand{\eemp}{\end{emphit}}

\newcommand{\bt}{\begin{tabular}}
\newcommand{\et}{\end{tabular}}

\newtheorem{myth}{Theorem}[section]
\newtheorem{mylem}{Lemma}[section]
\newtheorem{mycor}{Corollary}[section]

\newtheorem{mydef}{Definition}[section]
\newtheorem{myrem}{Remark}[section]

\begin{document}
\date{On the 50$^\textrm{th}$ Anniversary of Obafemi Awolowo University}
\title{Some simplicial complexes of universal Osborn loops \footnote{2010 Mathematics Subject
Classification. Primary 20N02, 20N05}
\thanks{{\bf Keywords and Phrases :} Osborn loops, universality, left universality, right
universality}}
\author{Ja\'iy\'e\d ol\'a T\`em\'it\d{\'o}p\d{\'e} Gb\d{\'o}l\'ah\`an \thanks{All correspondence to be addressed to this author.} \\
Department of Mathematics,\\
Obafemi Awolowo University, Ile Ife, Nigeria.\\
jaiyeolatemitope@yahoo.com, tjayeola@oauife.edu.ng}\maketitle
\begin{abstract}
A loop is shown to be a universal Osborn loop if and only if it has a particular simplicial complex. A loop is shown to be a universal Osborn loop and obeys two new identities if and only if it has another particular simplicial complex. A universal Osborn loop and four of its isotopes are shown to form a rectangular pyramid in a $3$-dimensional space.
\end{abstract}
\section{Introduction and Preliminaries}
A loop is called an Osborn loop if it obeys any of the two
identities below.
\begin{equation}\label{eq:4.1}
\textrm{OS$_3$}~:~(x\cdot yz)x=xy\cdot [(x^\lambda \cdot xz)\cdot x]
\end{equation}
\begin{equation}\label{eq:6.1}
\textrm{OS$_5$}~:~(x\cdot yz)x=xy\cdot [(x\cdot x^{\rho}z)\cdot x]
\end{equation}
For a comprehensive introduction to Osborn loops and its
universality, and a detailed literature review on it, readers should
check Jaiy\'e\d ol\'a , Ad\'en\'iran and S\`{o}l\'{a}r\`{i}n \cite{phd195} and Jaiy\'e\d ol\'a \cite{phd195b}.
In this present paper, we shall follow the style and notations used
in Jaiy\'e\d ol\'a , Ad\'en\'iran and S\`{o}l\'{a}r\`{i}n \cite{phd195} and Jaiy\'e\d ol\'a \cite{phd195b}. The only concepts
and notions which will be introduced here are those that were not
defined in Jaiy\'e\d ol\'a , Ad\'en\'iran and S\`{o}l\'{a}r\`{i}n \cite{phd195} and Jaiy\'e\d ol\'a \cite{phd195b}.
\begin{mydef}\label{definition:bijection}

Let $(L, \cdot )$ be a loop and $U, V, W\in SYM(L, \cdot )$.
\begin{enumerate}
\item If $(U, V, W)\in AUT(L, \cdot )$ for some $V, W$, then $U$ is called autotopic.
\item If $(U, V, W)\in AUT(L, \cdot )$ such that $W=U, V=I$, then $U$ is called $\lambda $-regular.
\item If $(U, V, W)\in AUT(L, \cdot )$ such that $U=I, W=V$, then $V$ is called $\rho $-regular.
\end{enumerate}
\end{mydef}
\paragraph{}
Drisko \cite{dris} while considering the action of isotopisms and autotopisms of loops, found it convenient to think of a loop $\mathcal{Q}=(Q, \cdot ,\backslash ,/)$ in terms of the set $T_Q$ of all ordered triples $(x,y,z)$ of elements of $Q$ such that $x\cdot y=z$. An isotopism $(\alpha ,\beta ,\gamma )$ from $G$ to $H$ takes $(x,y,z)\in T_G$ to $(x\alpha,y\beta,z\gamma)\in T_H$. We shall adopt his conventions at some points in time. We shall denote by $[\alpha,\beta]$, the commutator of any $\alpha,\beta\in SYM(G,\cdot )$.

Let $(Q, \cdot ,\backslash ,/)$ be a loop, then we shall be making use of the following notations for principal isotopes of $(Q, \cdot)$.
\begin{itemize}
\item $(Q,\ast_0 )$ represents $Q_{x,v}$;
\item $(Q,\circ_0 )$ represents $Q_{_{u,\phi_0(x,u,v)}}$, $\phi_0(x,u,v)=(u\backslash
([(uv)/(u\backslash (xv))]v))$;
\item $(Q,\circ_1 )$ represents $Q_{_{u,[u\backslash (xv)]}}$;
\item $(Q,\ast_1 )$ represents $Q_{_{\phi_1(x,u,v),v}}$, $\phi_1(x,u,v)=(u\backslash
([(uv)/(u\backslash (xv))]v))$ for all
$x,u,v\in Q$;
\item $(Q,\circ_2)$ represents $Q_{_{x,\phi_2(x,u,v)}}$, $\phi_2(x,u,v)=(u\backslash[(u/v)(u\backslash
(xv))])$;
\item $(Q,\circ_3 )$ represents $Q_{_{[x\cdot
u\backslash v]/v,[u\backslash (xv)]}}$;
\item $(Q,\ast_2)$ represents $Q_{u,e}$;
\item $(Q,\ast_3)$ represents $Q_{e,v}$.
\end{itemize}
Let $(G,\cdot )$ be a loop and let
$$BS_2(G,\cdot )=\{\theta\in SYM(G)~:~G(a,b)\overset{\theta}{\cong}G(c,d)~\textrm{for some}~a,b,c,d\in
G\}.$$ As shown in Bryant and Schneider \cite{phd92},
$BS_2(G,\cdot)$ forms a group for a loop $(G,\cdot )$ and it shall
be called the second Bryant-Schneider group (2$^{\textrm{nd}}$ BSG)
of the loop.

Consider the following two notions in algebraic topology.
\begin{mydef}
Let $V_Q$ be a set of isotopes of a loop $(Q,\cdot )$ and let $S_Q\subseteq ${\large
2}$^{{}^{V_Q}}$ such that $\phi\in S_Q$. If $S_Q$ is a topology on $V_Q$, then it is called the topology of isotopes of the loop $Q$ and the pair $(V_Q,S_Q)$ is called a topological space of isotopes of $Q$ if $(V_Q,S_Q)$ is a topological space.
\end{mydef}
Based on the above notion of topological space of isotopes of a loop, the following facts are direct consequences.
\begin{mylem}
Let $(Q,\cdot )$ be a loop and let $V_Q$ be the set of isotopes of $Q$.
Then, $\bigg(V_Q,${\large
2}$^{{}^{V_Q}}\bigg)$ is a topological space of isotopes of $Q$.
\end{mylem}
\begin{mylem}
Let $(Q,\cdot )$ be a G-loop and let $V_Q$ be the set of isotopes of $Q$. Let $S_Q=\{X_i\}_{i\in\Omega}\subseteq ${\large
2}$^{{}^{V_Q}}$ such that $\phi\in S_Q$ and $x_{i_j}\cong x_{i_k}$ for all $x_{i_j},x_{i_k}\in X_i$.
Then, $(V_Q,S_Q)$ is a topological space of isotopes of $Q$.
\end{mylem}
\begin{mycor}
Let $(Q,\cdot )$ be a CC-loop or VD-loop or K-loop or Buchsteiner loop or extra loop or group. Let $S_Q=\{X_i\}_{i\in\Omega}\subseteq ${\large
2}$^{{}^{V_Q}}$ such that $\phi\in S_Q$ and $x_{i_j}\cong x_{i_k}$ for all $x_{i_j},x_{i_k}\in X_i$.
Then, $(V_Q,S_Q)$ is a topological space of isotopes of $Q$.
\end{mycor}
\begin{mydef}
A simplicial complex is a pair $(V,S)$ where $V$ is a set of points called vertices and $S$ is a given family of finite subsets, called simplexes, so that the following  conditions are satisfied:
 \begin{enumerate}
 \item all points of $V$ are simplexes;
 \item any non-empty subset of a simplex is a simplex.
 \end{enumerate}
 A simplex consisting of $(n+1)$ points is called $n$-dimensional simplex.
\end{mydef}
\begin{mydef}
Let $V_Q$ be a set of isotopes of a loop $(Q,\cdot )$ and let $S_Q\subseteq ${\large
2}$^{{}^{V_Q}}$. If $K_Q=(V_Q,S_Q)$ is a simplicial complex, then $K_Q$ is called a trivial simplicial complex of isotopes of the loop $Q$.
\end{mydef}
\begin{mydef}
Let $V_Q$ be a set of isotopes of a loop $(Q,\cdot )$ and let $S_Q=\{X_i\}_{i\in\Omega}\subseteq ${\large
2}$^{{}^{V_Q}}$ such that $x_{i_j}\cong x_{i_k}$ for all $x_{i_j},x_{i_k}\in X_i$. If $K_Q=(V_Q,S_Q)$ is a simplicial complex, then $K_Q$ is called a non-trivial simplicial complex of isotopes or simplicial complex of isotopes of the loop $Q$.
\end{mydef}
The facts below follow suite.
\begin{mylem}
Let $(Q,\cdot )$ be a loop and let $V_Q$ be the set of isotopes of $Q$.
Then, $\bigg(V_Q,${\large
2}$^{{}^{V_Q}}\bigg)$ is a trivial simplicial complex of isotopes of $Q$.
\end{mylem}
\begin{mylem}
Let $(Q,\cdot )$ be a G-loop and let $V_Q$ be the set of isotopes of $Q$. Let $S_Q=\{X_i\}_{i\in\Omega}\subseteq ${\large
2}$^{{}^{V_Q}}$ such that $x_{i_j}\cong x_{i_k}$ for all $x_{i_j},x_{i_k}\in X_i$.
Then, $(V_Q,S_Q)$ is a simplicial complex of isotopes of $Q$.
\end{mylem}
\begin{mycor}
Let $(Q,\cdot )$ be a CC-loop or VD-loop or K-loop or Buchsteiner loop or extra loop or group. Let $S_Q=\{X_i\}_{i\in\Omega}\subseteq ${\large
2}$^{{}^{V_Q}}$ such that $x_{i_j}\cong x_{i_k}$ for all $x_{i_j},x_{i_k}\in X_i$.
Then, $(V_Q,S_Q)$ is a simplicial complex of isotopes of $Q$.
\end{mycor}

\begin{mydef}
Let $K=(V,S)$ and $K'=(V',S')$ be two simplicial complexes. A simplicial map $f~:~K\to K'$ is a set map $f~:~V\to V'$ satisfying the property: for every simplex $x\in S$, the image $f(x)\in S'$.
\end{mydef}

In this work, the notion of simplicial complex is used to characterize universal Osborn loops. The following results are important for the set objective.
\begin{myth}\label{1:4}(Jaiy\'e\d ol\'a , Ad\'en\'iran and S\`{o}l\'{a}r\`{i}n \cite{phd195})

Let $\mathcal{Q}=(Q, \cdot ,\backslash ,/)$ be a loop and $\gamma_0
(x,u,v)=\mathbb{R}_vR_{[u\backslash (xv)]}\mathbb{L}_uL_x$ for all
$x,u,v\in Q$, then $\mathcal{Q}$ is a universal Osborn
loop if and only if the commutative diagram
\begin{equation}\label{eq:7}
\begin{diagram}
&                                                                               &(Q,\circ_0 )\\
&\ruTo^{(R_{\phi_0 (x,u,v)},L_u,I)}_{}&\dTo^{(\gamma_0,\gamma_0,\gamma_0)}_{\textrm{isomorphism}}\\
(Q,\cdot ) &\rTo^{(R_v,L_x,I)}_{\textrm{principal isotopism}} &(Q,\ast_0)
\end{diagram}
\end{equation}
holds.
\end{myth}
\begin{myth}\label{post1:4}(Ja\'iy\'e\d ol\'a \cite{phd195b})

Let $\mathcal{Q}=(Q, \cdot ,\backslash ,/)$ be a loop and $\gamma_1
(x,u,v)=\mathbb{R}_vR_{[u\backslash (xv)]}\mathbb{L}_uL_x$ for all
$x,u,v\in Q$, then $\mathcal{Q}$ is a universal Osborn
loop if and only if the commutative diagram
\begin{equation}\label{eq:8}
\begin{diagram}
&                                                                               &(Q,\ast_1 )\\
&\ruTo^{(R_v,L_{\phi_1 (x,u,v)},I)}_{}&\dTo^{(\gamma_1,\gamma_1,\gamma_1)}_{\textrm{isomorphism}}\\
(Q,\cdot ) &\rTo^{(R_{[u\backslash (xv)]},L_u,I)}_{\textrm{principal isotopism}} &(Q,\circ_1)
\end{diagram}
\end{equation}
holds.
\end{myth}

\begin{myth}\label{1:12}(Jaiy\'e\d ol\'a , Ad\'en\'iran and S\`{o}l\'{a}r\`{i}n \cite{phd195})

Let $\mathcal{Q}=(Q, \cdot ,\backslash ,/)$ be a loop and $\gamma_0
(x,u,v)=\mathbb{R}_vR_{[u\backslash (xv)]}\mathbb{L}_uL_x$ for all
$x,u,v\in Q$, then $\mathcal{Q}$ is a universal Osborn
loop implies the commutative diagram
\begin{equation}\label{eq:7.m}
\begin{diagram}
&                                                                               &(Q,\circ_2 )\\
&\ruTo^{(R_{\phi_2 (x,u,v)},L_x,I)}_{}&\uTo^{(\gamma_0,\gamma_0,\gamma_0)}_{\textrm{isomorphism}}\\
(Q,\cdot ) &\rTo^{(I,L_u,I)}_{\textrm{principal isotopism}} &(Q,\ast_2)
\end{diagram}
\end{equation}
holds.
\end{myth}

\begin{myth}\label{post1:12}(Ja\'iy\'e\d ol\'a \cite{phd195b})

Let $\mathcal{Q}=(Q, \cdot ,\backslash ,/)$ be a loop and $\gamma_1
(x,u,v)=\mathbb{R}_vR_{[u\backslash (xv)]}\mathbb{L}_uL_x$ for all
$x,u,v\in Q$, then $\mathcal{Q}$ is a universal Osborn
loop implies the commutative diagram
\begin{equation}\label{eq:8.m}
\begin{diagram}
&                                                                               &(Q,\circ_3 )\\
&\ruTo^{(R_{[u\backslash (xv)]},L_{[x\cdot
u\backslash v]/v},I)}_{}&\uTo^{(\gamma_1,\gamma_1,\gamma_1)}_{\textrm{isomorphism}}\\
(Q,\cdot ) &\rTo^{(R_v,I,I)}_{\textrm{principal isotopism}} &(Q,\ast_3)
\end{diagram}
\end{equation}
holds.
\end{myth}

\begin{mylem}\label{drispost}(Drisko \cite{dris})

Let $\mathcal{Q}=(Q, \cdot ,\backslash ,/)$ be a loop. Then $Q_{f,g}\cong Q_{c,d}$ if and only if there exists $(\alpha , \beta ,\gamma )\in AUT(\mathcal{Q})$ such that $(f,g,fg)(\alpha ,\beta ,\gamma )=(c,d,cd)$.
\end{mylem}

\begin{myth}\label{0:3}(Bryant and Schneider \cite{phd92})

Let $(Q, \cdot ,\backslash ,/)$ be a quasigroup. If $Q_{a,b}\overset{I}{\cong}Q_{c,d}$ if and only if $c\cdot b,a\cdot d\in N_\mu (Q_{a,b})$ and $a\cdot b=c\cdot d$.
\end{myth}

\section{Main Results}
\begin{myth}\label{2post1:10}
Let $\mathcal{Q}=(Q, \cdot ,\backslash ,/)$ be a universal Osborn loop. Then, the following are necessary and sufficient for each other.
\begin{enumerate}
\item $(Q,\circ_0 )\overset{I}{\cong}(Q,\circ_1 )$.
\item $(Q,\ast_0 )\overset{I}{\cong}(Q,\ast_1 )$.
\item $\mathcal{Q}$ is a boolean group.
\end{enumerate}
\end{myth}
{\bf Proof}\\
By combining the commutative diagrams in Equation~\ref{eq:7} and Equation~\ref{eq:8}, we have the commutative diagram below.
\begin{equation}\label{eq:9}
\begin{diagram}
(Q,\circ_1)   & &  &     (Q,\circ_1)&&&&&&&&(Q,\circ_1)\\
    &&&&\\
    &&&&\\
    & &  & \uTo(0,2)^{\gamma_{01}^\circ}&\\
    & &  & (Q,\circ_0)\\
    \uTo^{(R_{[u\backslash (xv)]},L_u,I)}&\ruTo^{(R_{\phi_0},L_u,I)} &  &    &\rdTo(2,2)^{\gamma_0}&\\
   (Q,\cdot )&\rTo^{(R_v,L_x,I)}&&&&(Q,\ast_0)& &\\
\dTo^{(R_v,L_{\phi_1},I)}    &&&&&&\rdTo(6,3)^{\gamma_{01}^*}&\\
&  &    &   &&&&&&&&\uTo(7,0)_{\gamma_1}\\
(Q,\ast_1) &  &    & &  &&&&&&&(Q,\ast_1)&
\end{diagram}
\end{equation}
Let
\begin{displaymath}
(Q,\circ_0 )\xrightarrow[\textrm{isotopism}]{(\delta_{01}^\circ ,\varepsilon_{01}^\circ ,\pi_{01}^\circ)}(Q,\circ_1 ).
\end{displaymath}
So, from Equation~\ref{eq:9},
\begin{gather*}
(R_{\phi_0 (x,u,v)},L_u,I)(\delta_{01}^\circ ,\varepsilon_{01}^\circ ,\pi_{01}^\circ)=(R_{[u\backslash (xv)]},L_u,I)\Rightarrow\\
(R_{\phi_0 (x,u,v)}\delta_{01}^\circ ,L_u\varepsilon_{01}^\circ,\pi_{01}^\circ)=(R_{[u\backslash (xv)]},L_u,I)\Leftrightarrow\\
R_{\phi_0 (x,u,v)}\delta_{01}^\circ =R_{[u\backslash (xv)]},~L_u\varepsilon_{01}^\circ =L_u~\textrm{and}~\pi_{01}^\circ =I\Leftrightarrow\\
\delta_{01}^\circ =R_{\phi_0 (x,u,v)}^{-1}R_{[u\backslash (xv)]},~\varepsilon_{01}^\circ =L_u^{-1}L_u=I~\textrm{and}~\pi_{01}^\circ =I.
\end{gather*}
Thus, $(Q,\circ_0 )\cong (Q,\circ_1 )$ iff $\delta_{01}^\circ =\varepsilon_{01}^\circ =I$ iff
\begin{gather*}
R_{\phi_0 (x,u,v)}^{-1}R_{[u\backslash (xv)]}=I\Leftrightarrow \phi_0 (x,u,v)=[u\backslash (xv)]\\
(u\backslash
([(uv)/(u\backslash (xv))]v))=[u\backslash (xv)]\Leftrightarrow x\backslash (uv)=u\backslash (xv).
\end{gather*}
Similarly, by using the procedure above, it can be shown that $(Q,\ast_0 )\cong (Q,\ast_1 )$ iff $x\backslash (uv)=u\backslash (xv)$.

Keeping in mind that every Osborn loop of exponent 2 is an abelian group, hence, a Boolean group. This completes the proof.
\begin{myrem}
It can be observed that in a universal Osborn loop $\mathcal{Q}=(Q, \cdot ,\backslash ,/)$ and for $\gamma_0
(x,u,v)$ and $\gamma_1(x,u,v)$ of Theorem~\ref{1:4} and Theorem~\ref{post1:4},
$\gamma_0
(x,u,v)=\gamma_1(x,u,v)$ if and only if $\Big[\mathbb{L}_uL_x,\mathbb{R}_vR_{[u\backslash (xv)]}\Big]=I$ for all $x,u,v\in Q$.

The proof of Theorem~\ref{2post1:10} can also be achieved by making use of Theorem~\ref{0:3}. Take $a=u$, $b=\phi_0 (x,u,v)$, $c=u$ and $d=u\backslash (xv)$. Then, $(Q,\circ_0 )\overset{I}{\cong}(Q,\circ_1 )$ iff
\begin{gather*}
\textrm{(i)}~u\phi_0 (x,u,v)\in N_\mu \big((Q,\circ_0 )\big),~\textrm{(ii)}~u[u\backslash (xv)]\in N_\mu \big((Q,\circ_0 )\big),~\textrm{(iii)}~u\phi_0 (x,u,v)=u[u\backslash (xv)]\Leftrightarrow
\end{gather*}
$\mathcal{Q}$ is a Boolean group.
\end{myrem}

\begin{myth}\label{2post1:11}
Let $\mathcal{Q}=(Q, \cdot ,\backslash ,/)$ be a universal Osborn loop. Then
$(Q,\circ_0 )\cong(Q,\circ_1 )$ if and only if there exists $(I,\beta ,\gamma )\in AUT(\mathcal{Q})$ such that
\begin{equation}\label{eq:11}
uv=xR_v\mathbb{L}_u\beta^{-1}L_u\mathbb{R}_v\cdot xR_v\mathbb{L}_u=xR_v\gamma^{-1}\mathbb{R}_v\cdot xR_v\mathbb{L}_u
\end{equation}
for all $x,u,v\in Q$.
\end{myth}
{\bf Proof}\\
Following Lemma~\ref{drispost},
$(Q,\circ_0 )\cong(Q,\circ_1 )$ if and only if there exists $(\alpha,\beta ,\gamma )\in AUT(\mathcal{Q})$ such that
\begin{gather*}
(u,\phi_0 (x,u,v),u\phi_0 (x,u,v))(\alpha,\beta ,\gamma )=(u,[u\backslash (xv)],xv)\Leftrightarrow\\
(u\alpha,\phi_0 (x,u,v)\beta,(u\phi_0 (x,u,v))\gamma )=(u,[u\backslash (xv)],xv)\Leftrightarrow\\
u\alpha =u,~\phi_0 (x,u,v)\beta =[u\backslash (xv)]~\textrm{and}~(u\phi_0 (x,u,v))\gamma=xv\Leftrightarrow\\
\alpha =I,~\{u\backslash
([(uv)/(u\backslash (xv))]v)\}\beta =u\backslash (xv)~\textrm{and}~\{
[(uv)/(u\backslash (xv))]v\}\gamma=xv\Leftrightarrow\\
\alpha =I,~[(uv)/(u\backslash (xv))]R_v\mathbb{L}_u\beta =xR_v\mathbb{L}_u~\textrm{and}~[(uv)/(u\backslash (xv))]R_v\gamma=xR_v\Leftrightarrow\\
\alpha =I,~(uv)/(u\backslash (xv))=xR_v\mathbb{L}_u\beta^{-1}L_u\mathbb{R}_v
~\textrm{and}~[(uv)/(u\backslash (xv))]=xR_v\gamma^{-1}\mathbb{R}_v\Leftrightarrow\\
\alpha =I,~uv=xR_v\mathbb{L}_u\beta^{-1}L_u\mathbb{R}_v\cdot xR_v\mathbb{L}_u
~\textrm{and}~uv=xR_v\gamma^{-1}\mathbb{R}_v\cdot xR_v\mathbb{L}_u\Leftrightarrow
\end{gather*}
there exists $(I,\beta ,\gamma )\in AUT(\mathcal{Q})$ such that
\begin{displaymath}
uv=xR_v\mathbb{L}_u\beta^{-1}L_u\mathbb{R}_v\cdot xR_v\mathbb{L}_u=xR_v\gamma^{-1}\mathbb{R}_v\cdot xR_v\mathbb{L}_u.
\end{displaymath}

\begin{myrem}
If the autotopism $(\alpha,\beta ,\gamma )$ in Theorem~\ref{2post1:11} is the identity autotopism, then we shall have the equivalence of 1. and 3. of Theorem~\ref{2post1:10}.
\end{myrem}

\begin{mycor}\label{2post1:12}
Let $\mathcal{Q}=(Q, \cdot ,\backslash ,/)$ be a universal Osborn loop. Then
$(Q,\circ_0 )\cong(Q,\circ_1 )$ implies that there exists $(I,\beta ,\gamma )\in AUT(\mathcal{Q})$ such that $\gamma =\mathbb{L}_u\beta L_u$ for all $u\in Q$. Hence,
\begin{enumerate}
\item $\gamma =\beta$ iff $[\beta ,L_u]=I$ or $[\gamma ,L_u]=I$. Thence, $\beta$ is a $\rho$-regular permutation.
\item  $\gamma =L_u$ iff $\beta =L_u$. Thence, $\mathcal{Q}$ is an abelian group.
\end{enumerate}
\end{mycor}
{\bf Proof}\\
The proof of these follows from the fact in Theorem~\ref{2post1:11} that
\begin{displaymath}
xR_v\mathbb{L}_u\beta^{-1}L_u\mathbb{R}_v\cdot xR_v\mathbb{L}_u=xR_v\gamma^{-1}\mathbb{R}_v\cdot xR_v\mathbb{L}_u\Rightarrow
\end{displaymath}
$\mathbb{L}_u\beta L_u=\gamma$ for all $u\in Q$.

\begin{myth}\label{2post1:11b}
Let $\mathcal{Q}=(Q, \cdot ,\backslash ,/)$ be a universal Osborn loop. Then
$(Q,\ast_0 )\cong(Q,\ast_1 )$ if and only if there exists $(\delta,I ,\pi )\in AUT(\mathcal{Q})$ such that
\begin{equation}\label{eq:11b}
uv=x\cdot x\delta R_v\mathbb{L}_u=x\cdot xR_v\pi\mathbb{L}_u
\end{equation}
for all $x,u,v\in Q$.
\end{myth}
{\bf Proof}\\
Following Lemma~\ref{drispost},
$(Q,\ast_0 )\cong(Q,\ast_1 )$ if and only if there exists $(\delta,\varepsilon ,\pi )\in AUT(\mathcal{Q})$ such that $(x,v,xv)(\delta,\varepsilon ,\pi )=(\phi_1(x,u,v),v,\phi_1(x,u,v)v)$. The procedure of the proof of the remaining part is similar to that of Theorem~\ref{2post1:11}.

\begin{myrem}
If the autotopism $(\delta,\varepsilon ,\pi )$ in Theorem~\ref{2post1:11b} is the identity autotopism, then we shall have the equivalence of 2. and 3. of Theorem~\ref{2post1:10}.
\end{myrem}

\begin{mycor}\label{2post1:13}
Let $\mathcal{Q}=(Q, \cdot ,\backslash ,/)$ be a universal Osborn loop. Then
$(Q,\ast_0 )\cong(Q,\ast_1 )$ implies that there exists $(\delta,I ,\pi )\in AUT(\mathcal{Q})$ such that $\pi =\mathbb{R}_v\delta R_v$ for all $v\in Q$. Hence,
\begin{enumerate}
\item $\pi =\delta$ iff $[\delta ,R_v]=I$ or $[\pi ,R_v]=I$. Thence, $\delta$ is a $\lambda$-regular permutation.
\item  $\delta =R_v$ iff $\pi =R_v$. Thence, $\mathcal{Q}$ is an abelian group.
\end{enumerate}
\end{mycor}
{\bf Proof}\\
The proof of these follows from the fact in Theorem~\ref{2post1:11b} that
\begin{displaymath}
x\cdot x\delta R_v\mathbb{L}_u=x\cdot xR_v\pi\mathbb{L}_u\Rightarrow
\end{displaymath}
$\pi =\mathbb{R}_v\delta R_v$ for all $v\in Q$.

\begin{myth}\label{2post1:14}
Let $\mathcal{Q}=(Q, \cdot ,\backslash ,/)$ be a universal Osborn loop. Then
$(Q,\circ_0 )\cong(Q,\circ_1 )$ and $(Q,\ast_0 )\cong(Q,\ast_1 )$ if and only if there exists $(I,\beta ,\gamma ),(\delta,I ,\pi )\in AUT(\mathcal{Q})$ such that
\begin{equation}\label{eq:10}
uv=xR_v\mathbb{L}_u\beta^{-1}L_u\mathbb{R}_v\cdot xR_v\mathbb{L}_u=xR_v\gamma^{-1}\mathbb{R}_v\cdot xR_v\mathbb{L}_u=x\cdot x\delta R_v\mathbb{L}_u=x\cdot xR_v\pi\mathbb{L}_u
\end{equation}
for all $x,u,v\in Q$
\end{myth}
{\bf Proof}\\
This is achieved by simply combining Theorem~\ref{2post1:11} and Theorem~\ref{2post1:11b}.

\begin{myth}\label{2post1:15}
Let $\mathcal{Q}=(Q, \cdot ,\backslash ,/)$ be a universal Osborn loop. If
$(Q,\circ_0 )\overset{\gamma_{01}^\circ}{\cong}(Q,\circ_1 )$ and $(Q,\ast_0 )\overset{\gamma_{01}^\ast}{\cong}(Q,\ast_1 )$, then $\gamma_0\gamma_{01}^\ast\gamma_1=\gamma_{01}^\circ$.
\end{myth}
{\bf Proof}\\
The commutative diagram in Equation~\ref{eq:9} proves this.

\begin{mycor}\label{2post1:16}
Let $\mathcal{Q}=(Q, \cdot ,\backslash ,/)$ be a universal Osborn loop. If $(Q,\circ_0 )\cong(Q,\circ_1 )$ and $(Q,\ast_0 )\cong(Q,\ast_1 )$, then the following are necessary and sufficient for each other.
\begin{multicols}{3}
\begin{enumerate}
\item $\beta =I$.
\item $\gamma =I$.
\item $\delta =I$.
\item $\pi =I$.
\item $(Q,\circ_0 )\overset{I}{\cong}(Q,\circ_1 )$.
\item $(Q,\ast_0 )\overset{I}{\cong}(Q,\ast_1 )$.
\item $\mathcal{Q}$ is a boolean group.
\end{enumerate}
\end{multicols}
\end{mycor}
{\bf Proof}\\
To prove the equivalence of 1. to 4. and 7., use Equation~\ref{eq:10} of
Theorem~\ref{2post1:14}. The proof of the equivalence of 5. to 7. follows from Theorem~\ref{2post1:10}.

\begin{myrem}
Corollary~\ref{2post1:16} is a very important result in this study. It gives us the main distinctions between Theorem~\ref{2post1:10} and Theorem~\ref{2post1:14}. That is, the
necessary and sufficient condition(s) under which the isomorphisms $(Q,\circ_0 )\cong(Q,\circ_1 )$ and $(Q,\ast_0 )\cong(Q,\ast_1 )$ will be trivial. And the condition(s) is when any of the autotopic permutations of $\beta$, $\gamma$, $\delta$ and
$\pi$ of Theorem~\ref{2post1:11} and Theorem~\ref{2post1:11b} is equal to the identity mapping.
\end{myrem}
\paragraph{}
Next, it is important to deduce the actual definitions of the autotopic mappings $\beta$, $\gamma$, $\delta$, $\pi$ and the isomorphisms $\gamma_{01}^\ast$ and $\gamma_{01}^\circ$. Recall that by the necessary part of Lemma~\ref{drispost}, if $\mathcal{Q}=(Q, \cdot ,\backslash ,/)$ is a loop and $Q_{f,g}\overset{\theta}{\cong}Q_{c,d}$, then there exists $(A,B,C)\in AUT(\mathcal{Q})$ such that $(f,g,fg)(A,B,C)=(c,d,cd)$. According to the proof of this,
\begin{equation}\label{eq:10.m}
(A,B,C)=(R_g\theta R_d^{-1},L_f\theta L_c^{-1},\theta)\Leftrightarrow A=R_g\theta R_d^{-1},~B=L_f\theta L_c^{-1}~\textrm{and}~C=\theta.
\end{equation}
Thus,
\begin{gather*}
I=\alpha =R_{\phi_o(x,u,v)}\gamma_{01}^\circ R_{[u\backslash (xv)]}^{-1},~\beta=L_u\gamma_{01}^\circ L_u^{-1}~\textrm{and}~\gamma=\gamma_{01}^\circ\\
\gamma_{01}^\circ =\mathbb{R}_{\phi_o(x,u,v)}R_{[u\backslash (xv)]},~\beta=L_u\mathbb{R}_{\phi_o(x,u,v)}R_{[u\backslash (xv)]}\mathbb{L}_u^{-1}~\textrm{and}~\gamma=\mathbb{R}_{\phi_o(x,u,v)}R_{[u\backslash (xv)]}
\end{gather*}
and
\begin{gather*}
\delta =R_v\gamma_{01}^\ast R_v^{-1},~I=\varepsilon=L_x\gamma_{01}^\ast L_{\phi_1(x,u,v)}^{-1}~\textrm{and}~\pi=\gamma_{01}^\ast\\
\delta =R_v\gamma_{01}^\ast \mathbb{R}_v^{-1},~\gamma_{01}^\ast =\mathbb{L}_xL_{\phi_1(x,u,v)}~\textrm{and}~\pi=\gamma_{01}^\ast\\
\delta =R_v\mathbb{L}_xL_{\phi_1(x,u,v)}\mathbb{R}_v^{-1},~\gamma_{01}^\ast =\mathbb{L}_xL_{\phi_1(x,u,v)}~\textrm{and}~\pi=\mathbb{L}_xL_{\phi_1(x,u,v)}.
\end{gather*}
Therefore, Theorem~\ref{2post1:11} and Theorem~\ref{2post1:11b} can now be restated as follows.

\begin{myth}\label{2post1:17}
Let $\mathcal{Q}=(Q, \cdot ,\backslash ,/)$ be a universal Osborn loop. Then
$(Q,\circ_0 )\overset{\gamma_{01}^\circ}{\cong}(Q,\circ_1 )$ if and only if
\begin{equation}\label{eq:12}
y\cdot u\backslash [(uz)\psi_0]=(yz)\psi_0~\textrm{and}~uv=xR_v(R_v\psi_0)^{-1}\cdot xR_v\mathbb{L}_u
\end{equation}
where $\psi_0=\mathbb{R}_{\phi_o(x,u,v)}R_{[u\backslash (xv)]}$ for all $x,y,z,u,v\in Q$
\end{myth}
{\bf Proof}\\
Simply substitute
\begin{displaymath}
\beta=L_u\mathbb{R}_{\phi_o(x,u,v)}R_{[u\backslash (xv)]}\mathbb{L}_u^{-1}~\textrm{and}~\gamma=\mathbb{R}_{\phi_o(x,u,v)}R_{[u\backslash (xv)]}
\end{displaymath}
into Equation~\ref{eq:11} of Theorem~\ref{2post1:11}.

\begin{myth}\label{2post1:17b}
Let $\mathcal{Q}=(Q, \cdot ,\backslash ,/)$ be a universal Osborn loop. Then
$(Q,\ast_0 )\overset{\gamma_{01}^\ast}{\cong}(Q,\ast_1 )$ if and only if
\begin{equation}\label{eq:12b}
[(yv)\psi_1]/v\cdot z=(yz)\psi_1~\textrm{and}~uv=x\cdot u\backslash [(xv)\psi_1]
\end{equation}
where $\psi_1=\mathbb{L}_xL_{\phi_1(x,u,v)}$ for all $x,y,z,u,v\in Q$
\end{myth}
{\bf Proof}\\
Simply substitute
\begin{displaymath}
\delta =R_v\mathbb{L}_xL_{\phi_1(x,u,v)}\mathbb{R}_v^{-1}~\textrm{and}~\pi=\mathbb{L}_xL_{\phi_1(x,u,v)}
\end{displaymath}
into Equation~\ref{eq:11b} of Theorem~\ref{2post1:11b}.

\begin{mylem}\label{2post1:17c}
Let $\mathcal{Q}=(Q, \cdot ,\backslash ,/)$ be a loop.
\begin{enumerate}
\item $\mathcal{Q}$ is a universal Osborn loop and obeys Equation~\ref{eq:12} if and only if
$\gamma_0,\gamma_{01}^\circ\in BS_2(\mathcal{Q})$.
\item $\mathcal{Q}$ is a universal Osborn loop and obeys Equation~\ref{eq:12b} if and only if
$\gamma_1,\gamma_{01}^\ast\in BS_2(\mathcal{Q})$.
\end{enumerate}
\end{mylem}
{\bf Proof}\\
This follows by combining Theorem~\ref{1:4}, Theorem~\ref{post1:4}, Theorem~\ref{2post1:11} and Theorem~\ref{2post1:11b}

\begin{myrem}
It is a self exercise to confirm if $(Q,\circ_0 )\overset{\gamma_{01}^\circ}{\cong}(Q,\circ_1 )$ and $(Q,\ast_0 )\overset{\gamma_{01}^\ast}{\cong}(Q,\ast_1 )$ in some universal Osborn loops like Moufang loops and extra loops by simply verifying Equation~\ref{eq:12} and Equation~\ref{eq:12b}. Furthermore, the relation $\gamma_0\gamma_{01}^\ast\gamma_1=\gamma_{01}^\circ$ of Theorem~\ref{2post1:15} is justifiable as well. It must be noted also, that in any universal Osborn loop $\mathcal{Q}$, Equation~\ref{eq:12} and Equation~\ref{eq:12b} are necessary and sufficient conditions for $\gamma_{01}^\ast,\gamma_{01}^\circ\in BS_2(\mathcal{Q})$.
\end{myrem}

\paragraph{}
By combining the commutative diagrams in Equation~\ref{eq:7.m} and Equation~\ref{eq:8.m}, we have the commutative diagram below.
\begin{equation}\label{eq:17}
\begin{diagram}
(Q,\circ_3)   & &  &     (Q,\circ_3)&&&&&&&&(Q,\circ_3)\\
    &&&&\\
    &&&&\\
    & &  & \uTo(0,2)^{\gamma_{23}^\circ}&\\
    & &  & (Q,\circ_2)\\
    \uTo^{(R_{[u\backslash (xv)]},L_{\{[x\cdot u\backslash v]/v\}},I)}&\ruTo^{(R_{\phi_2},L_x,I)} &  &    &\rdTo(2,2)^{\gamma_0}&\\
   (Q,\cdot )&\rTo^{(I,L_u,I)}&&&&(Q,\ast_2)& &\\
\dTo^{(R_v,I,I)}    &&&&&&\rdTo(6,3)^{\gamma_{23}^*}&\\
&  &    &   &&&&&&&&\uTo(7,0)_{\gamma_1}\\
(Q,\ast_3) &  &    & &  &&&&&&&(Q,\ast_3)&
\end{diagram}
\end{equation}

\begin{myth}\label{2post1:18}
Let $\mathcal{Q}=(Q, \cdot ,\backslash ,/)$ be a universal Osborn loop. Then
$(Q,\circ_2 )\cong(Q,\circ_3 )$ if and only if there exists $(\lambda,\mu,\nu )\in AUT(\mathcal{Q})$ such that
\begin{equation}\label{eq:13}
\lambda =R_{u\backslash v}\mathbb{R}_v,~\mu=L_u\mathbb{L}_{u\backslash v}~\textrm{and}~
[x\cdot xR_v\mathbb{L}_u\mu^{-1}]\nu=x\lambda\cdot xR_v\mathbb{L}_u
\end{equation}
for all $x,u,v\in Q$.
\end{myth}
{\bf Proof}\\
Following Lemma~\ref{drispost},
$(Q,\circ_2 )\cong(Q,\circ_3 )$ if and only if there exists $(\lambda,\mu,\nu )\in AUT(\mathcal{Q})$ such that
$(x,\phi_2(x,u,v),x\phi_2(x,u,v))(\lambda,\mu,\nu )=([x\cdot
u\backslash v]/v,[u\backslash (xv)],\{[x\cdot
u\backslash v]/v\}[u\backslash (xv)])$. The procedure of the proof of the remaining part is similar to that of Theorem~\ref{2post1:11}.

\begin{mylem}\label{2post1:19}
Let $\mathcal{Q}=(Q, \cdot ,\backslash ,/)$ be a universal Osborn loop. Then
$(Q,\circ_2 )\overset{\gamma_{23}^\circ}{\cong}(Q,\circ_3 )$ if and only if there exists $(\lambda,\mu,\gamma_{23}^\circ )\in AUT(\mathcal{Q})$ such that
\begin{equation}\label{eq:14}
\gamma_{23}^\circ=\mathbb{R}_{\phi_2(x,u,v)}R_{u\backslash v}\mathbb{R}_vR_{[u\backslash (xv)]}=\mathbb{L}_xL_u\mathbb{L}_{u\backslash v}L_{\{[x\cdot
u\backslash v]/v\}}~\textrm{and}~
[x\cdot xR_v\mathbb{L}_u\mu^{-1}]\gamma_{23}^\circ=x\lambda\cdot xR_v\mathbb{L}_u
\end{equation}
for all $x,u,v\in Q$.
\end{mylem}
{\bf Proof}\\
Considering the commutative diagram in Equation~\ref{eq:17} and using Equation~\ref{eq:10.m},
\begin{displaymath}
\lambda=R_{\phi_2(x,u,v)}\gamma_{23}^\circ R_{[u\backslash (xv)]}^{-1},~\mu=L_x\gamma_{23}^\circ L_{\{[x\cdot
u\backslash v]/v\}}^{-1}~\textrm{and}~\nu=\gamma_{23}^\circ.
\end{displaymath}
The final conclusion follows from Theorem~\ref{2post1:18}.

\begin{mycor}\label{2post1:20}
Let $\mathcal{Q}=(Q, \cdot ,\backslash ,/)$ be a universal Osborn loop. $\gamma_{23}^\circ\in BS_2(\mathcal{Q})$ if and only if there exists $(\lambda,\mu,\gamma_{23}^\circ )\in AUT(\mathcal{Q})$ such that
\begin{equation}\label{eq:15}
\gamma_{23}^\circ=\mathbb{R}_{\phi_2(x,u,v)}R_{u\backslash v}\mathbb{R}_vR_{[u\backslash (xv)]}=\mathbb{L}_xL_u\mathbb{L}_{u\backslash v}L_{\{[x\cdot
u\backslash v]/v\}}~\textrm{and}~
[x\cdot xR_v\mathbb{L}_u\mu^{-1}]\gamma_{23}^\circ=x\lambda\cdot xR_v\mathbb{L}_u
\end{equation}
for all $x,u,v\in Q$.
\end{mycor}
{\bf Proof}\\
This follows from Lemma\ref{2post1:19}.

\begin{mycor}\label{2post1:21}
Let $\mathcal{Q}=(Q, \cdot ,\backslash ,/)$ be a loop. $\mathcal{Q}$ is a universal Osborn loop and $\gamma_{23}^\circ\in BS_2(\mathcal{Q})$ implies $\gamma_0\in BS_2(\mathcal{Q})$ and there exists $(\lambda,\mu,\gamma_{23}^\circ )\in AUT(\mathcal{Q})$ such that
\begin{equation}\label{eq:16}
\gamma_{23}^\circ=\mathbb{R}_{\phi_2(x,u,v)}R_{u\backslash v}\mathbb{R}_vR_{[u\backslash (xv)]}=\mathbb{L}_xL_u\mathbb{L}_{u\backslash v}L_{\{[x\cdot
u\backslash v]/v\}}~\textrm{and}~
[x\cdot xR_v\mathbb{L}_u\mu^{-1}]\gamma_{23}^\circ=x\lambda\cdot xR_v\mathbb{L}_u
\end{equation}
for all $x,u,v\in Q$.
\end{mycor}
{\bf Proof}\\
This follows from Theorem~\ref{1:12} and Lemma\ref{2post1:19}.

\subsection*{Simplicial Complex of Isotopes of a Universal Osborn Loop}
\begin{myth}\label{0}
Let $(Q,\cdot )$ be a loop. Let $V_0(Q)=\big\{(Q,\cdot ),(Q,\circ_0),(Q,\ast_0)\big\}$ and $S_0(Q)=\Big\{\{(Q,\cdot )\},\{(Q,\circ_0)\},\{(Q,\ast_0)\},\big\{(Q,\circ_0),(Q,\ast_0)\big\}\Big\}$. Then, $(Q,\cdot )$ is a universal Osborn loop if and only if $K_0(Q)=\Big(V_0(Q),S_0(Q)\Big)$ is a simplicial complex of isotopes of $(Q,\cdot )$.
\end{myth}
{\bf Proof}\\
This is proved with the help of Theorem~\ref{1:4}.

\begin{myth}\label{1}
Let $(Q,\cdot )$ be a loop. Let $V_1(Q)=\big\{(Q,\cdot ),(Q,\circ_1),(Q,\ast_1)\big\}$ and $S_1(Q)=\Big\{\{(Q,\cdot )\},\{(Q,\circ_1)\},\{(Q,\ast_1)\},\big\{(Q,\circ_1),(Q,\ast_1)\big\}\Big\}$. Then, $(Q,\cdot )$ is a universal Osborn loop if and only if $K_1(Q)=\Big(V_1(Q),S_1(Q)\Big)$ is a simplicial complex of isotopes of $(Q,\cdot )$.
\end{myth}
{\bf Proof}\\
This is proved with the help of Theorem~\ref{post1:4}.

\begin{myth}\label{2}
Let $(Q,\cdot )$ be a loop. Let $V_2(Q)=\big\{(Q,\cdot ),(Q,\circ_2),(Q,\ast_2)\big\}$ and $S_2(Q)=\Big\{\{(Q,\cdot )\},\{(Q,\circ_2)\},\{(Q,\ast_2)\},\big\{(Q,\circ_2),(Q,\ast_2)\big\}\Big\}$. If $(Q,\cdot )$ is a universal Osborn loop, then $K_2(Q)=\Big(V_2(Q),S_2(Q)\Big)$ is a simplicial complex of isotopes of $(Q,\cdot )$.
\end{myth}
{\bf Proof}\\
This is proved with Theorem~\ref{1:12}.

\begin{myth}\label{3}
Let $(Q,\cdot )$ be a loop. Let $V_3(Q)=\big\{(Q,\cdot ),(Q,\circ_3),(Q,\ast_3)\big\}$ and $S_3(Q)=\Big\{\{(Q,\cdot )\},\{(Q,\circ_3)\},\{(Q,\ast_3)\},\big\{(Q,\circ_3),(Q,\ast_3)\big\}\Big\}$. If $(Q,\cdot )$ is a universal Osborn loop, then $K_3(Q)=\Big(V_3(Q),S_3(Q)\Big)$ is a simplicial complex of isotopes of $(Q,\cdot )$.
\end{myth}
{\bf Proof}\\
This is proved with the aid of Theorem~\ref{post1:12}.

\begin{mycor}\label{01}
Let $(Q,\cdot )$ be a loop. Let $V_i(Q)=\big\{(Q,\cdot ),(Q,\circ_i),(Q,\ast_i)\big\}$ and $S_i(Q)=\Big\{\{(Q,\cdot )\},\{(Q,\circ_i)\},\{(Q,\ast_i)\},\big\{(Q,\circ_i),(Q,\ast_i)\big\}\Big\}$ for $i=0,1$. Then, $(Q,\cdot )$ is a universal Osborn loop if and only if $K_{01}(Q)=K_0(Q)\cup K_1(Q)=\bigg(V_0(Q)\cup V_1(Q),S_0(Q)\cup S_1(Q)\bigg)$ is a simplicial complex of isotopes of $(Q,\cdot )$.
\end{mycor}
{\bf Proof}\\
This follows from Theorem~\ref{0} and Theorem~\ref{1}.

\begin{mycor}\label{23}
Let $(Q,\cdot )$ be a loop. Let $V_i(Q)=\big\{(Q,\cdot ),(Q,\circ_i),(Q,\ast_i)\big\}$ and $S_i(Q)=\Big\{\{(Q,\cdot )\},\{(Q,\circ_i)\},\{(Q,\ast_i)\},\big\{(Q,\circ_i),(Q,\ast_i)\big\}\Big\}$ for $i=2,3$. If $(Q,\cdot )$ is a universal Osborn loop, then $K_{23}(Q)=K_2(Q)\cup K_3(Q)=\bigg(V_2(Q)\cup V_3(Q),S_2(Q)\cup S_3(Q)\bigg)$ is a simplicial complex of isotopes of $(Q,\cdot )$.
\end{mycor}
{\bf Proof}\\
This follows from Theorem~\ref{2} and Theorem~\ref{3}.

\begin{mycor}
Let $(Q,\cdot )$ be a loop. Let $V_i(Q)=\big\{(Q,\cdot ),(Q,\circ_i),(Q,\ast_i)\big\}$ and $S_i(Q)=\Big\{\{(Q,\cdot )\},\{(Q,\circ_i)\},\{(Q,\ast_i)\},\big\{(Q,\circ_i),(Q,\ast_i)\big\}\Big\}$ for $i=0,1,2,3$. If $(Q,\cdot )$ is a universal Osborn loop, then $\displaystyle K_{0123}(Q)=\bigcup^3_{i=0}K_i(Q)=\bigg(\bigcup^3_{i=0}V_i(Q),\bigcup^3_{i=0}S_i(Q)\bigg)$ is a simplicial complex of isotopes of $(Q,\cdot )$.
\end{mycor}
{\bf Proof}\\
This is proved by combining Corollary~\ref{01} and Corollary~\ref{23}.

\begin{myth}
Let $(Q,\cdot )$ be a loop. Let $V_{01}(Q)=\big\{(Q,\cdot ),(Q,\circ_0),(Q,\ast_0),(Q,\circ_1),(Q,\ast_1)\big\}$ and $S_{10}(Q)=\Big\{\{(Q,\cdot )\},\{(Q,\circ_0)\},\{(Q,\ast_0)\},\{(Q,\circ_1)\},\{(Q,\ast_1)\},\big\{(Q,\circ_0),(Q,\ast_0)\big\},\\
\big\{(Q,\circ_1),(Q,\ast_1)\big\},\big\{(Q,\circ_0),(Q,\circ_1)\big\},\big\{(Q,\ast_0),(Q,\ast_1)\big\},
\big\{(Q,\circ_0),(Q,\ast_1)\big\},\big\{(Q,\circ_1),(Q,\ast_0)\big\},\\
\big\{(Q,\circ_0),(Q,\circ_1),(Q,\ast_0)\big\},
\big\{(Q,\circ_0),(Q,\circ_1),(Q,\ast_1)\big\},\big\{(Q,\ast_0),(Q,\ast_1),(Q,\circ_0)\big\},\\
\big\{(Q,\ast_0),(Q,\ast_1),(Q,\circ_1)\big\},
\big\{(Q,\circ_0),(Q,\circ_1),(Q,\ast_0),(Q,\ast_1)\big\}\Big\}$. Then, $(Q,\cdot )$ is a universal Osborn loop and obey Equation~\ref{eq:12} and Equation~\ref{eq:12b} if and only if $K_{10}(Q)=\Big(V_{01}(Q),S_{10}(Q)\Big)$ is a simplicial complex of isotopes of $(Q,\cdot )$.
\end{myth}
{\bf Proof}\\
This is proved with the aid of Theorem~\ref{0}, Theorem~\ref{1}, Theorem~\ref{2post1:17} and Theorem~\ref{2post1:17b}.

\begin{myth}\label{3.1}
Let $(Q,\cdot )$ be a universal Osborn loop. Let $V_i(Q)=\big\{(Q,\cdot ),(Q,\circ_i),(Q,\ast_i)\big\}$, $S_i(Q)=\Big\{\{(Q,\cdot )\},\{(Q,\circ_i)\},\{(Q,\ast_i)\},\big\{(Q,\circ_i),(Q,\ast_i)\big\}\Big\}$ and $K_i=(V_i(Q),S_i(Q))$ for $i=0,1,2,3$.
Define $f_{ij}~:~K_i\to K_j$ as \begin{equation*}
f_{ij}~:~\begin{cases}
(Q,\cdot )  & \longmapsto  (Q,\cdot )\\
 (Q,\circ_i) &  \longmapsto (Q,\circ_j)\\
 (Q,\ast_i) & \longmapsto  (Q,\ast_j)
\end{cases}
~~i,j=0,1,2,3~\textrm{such that}~i\ne j.
\end{equation*}
Then, $f_{ij}$ is a simplicial map.
\end{myth}
{\bf Proof}\\
This is proved by Theorem~\ref{0}, Theorem~\ref{1}, Theorem~\ref{2} and Theorem~\ref{3}.

\begin{myth}
Let $(G,\cdot )$ and $(H,\star )$ be two loop isotopes under the triple $(A,B,C)$. For $D\in\{A,B,C\}$, if $D=E_1E_2\cdots E_i\cdots E_n$,~$E_i~:~G\to H,~i=1,\cdots n$ been bijections such that there does not exist $r\ge n$ for which $D=E_1E_2\cdots E_i\cdots E_r$, then the length of $D$, $|D|=n~\textrm{units}$. If $D=I$, the identity mapping, then $|D|=0$. The length of the isotopism $(G,\cdot
)\xrightarrow[\textrm{Isotopism}]{(A,B,C)}(H,\star )$ is giving by $|(A,B,C)|=|A|+|B|+|C|~\textrm{units}$.
For an isotopism $(G,\cdot
)\xrightarrow[\textrm{Isotopism}]{(A,B,C)}(H,\star )$, let the two loops $(G,\cdot )$ and $(H,\star )$ represent points in a $3$-dimensional space and let an isotopism from $(G,\cdot )$ to $(H,\star )$ be a line with $(G,\cdot )$ and $(H,\star )$ as end-points. The set of loops $V_{01}(Q)=\big\{(Q,\cdot ),(Q,\circ_0),(Q,\ast_0),(Q,\circ_1),(Q,\ast_1)\big\}$ where $(Q,\cdot )$ is a universal Osborn loop, form a rectangular pyramid with apex $(Q,\cdot )$.
\end{myth}
{\bf Proof}\\
We shall make use of the combined commutative diagram (\ref{eq:9}) as shown in the proof of Theorem~\ref{2post1:10}. There are four isotopes of $(Q,\cdot )$ as shown in the combined commutative diagram (\ref{eq:9}), namely $(Q,\circ_i),(Q,\ast_i)$ for $i=0,1$. The length of each of the isotopisms $(R_{[u\backslash (xv)]},L_u,I),(R_{\phi_0},L_u,I),(R_v,L_{\phi_1},I),(R_v,L_x,I)$ is $2~\textrm{units}$. The length of each of the isomorphisms $\gamma_0
(x,u,v)=\mathbb{R}_vR_{[u\backslash (xv)]}\mathbb{L}_uL_x$ and $\gamma_1
(x,u,v)=\mathbb{R}_vR_{[u\backslash (xv)]}\mathbb{L}_uL_x$ is $12~\textrm{units}$. The length of each of the isomorphisms $\gamma_{01}^\circ =\mathbb{R}_{\phi_o(x,u,v)}R_{[u\backslash (xv)]}$ and $\gamma_{01}^\ast =\mathbb{L}_xL_{\phi_1(x,u,v)}$ is $6~\textrm{units}$. Hence, the four loop isotopes $(Q,\circ_i),(Q,\ast_i)$ for $i=0,1$ of $(Q,\cdot )$ form a rectangle. Thus, taking $(Q,\cdot )$ as an apex and the four isotopism as lines drawn from the apex to the four vertices of the rectangle, we have a rectangular pyramid.

\end{document}